# Number of Edges with Shortest Cycle $k$ in a Kautz Graph

Vance Faber

Hoquiam, Washington

**Abstract**

For the Kautz digraph $K(d, D)$, let $\rho_k(d, D)$ be the number of oriented edges whose shortest directed cycle has length $k+1$, and define $\Delta_k(d, D) := \rho_k(d, D) - \rho_k(d, D-1)$. We give an exact, finite-dimensional matrix product that computes $\Delta_k(d, D)$ *directly*, without first computing $\rho$. In particular,
$$\Delta_k(d, D) = 0 \quad \text{for } k < \lfloor D/2 \rfloor + 2.$$
and $\Delta_k(d, D)$ is positive for every $k$ with $\lfloor D/2 \rfloor + 2 \leq k \leq D - 1$.

## 1 Background and the $\sigma/\rho$ identity

Let $K(d, D)$ be the Kautz digraph of out-degree $d$ and diameter $D$. For distinct vertices $u, v$, let $\text{dist}(u, v)$ denote the directed distance and let $Z_{uv}$ be the unique shortest directed cycle containing the edge $uv$. Define, for $1 \leq k \leq D$,

$$\sigma_k := \#\{(u, v) : \text{dist}(u, v) = k\}, \qquad \rho_k := \#\{uv \in E : Z_{uv} \text{ has length } k+1\}.$$

The following counting identity is folklore; we include a short proof for completeness. See Aigner [2] for fundamentals of directed line graphs and Fiol–Yebra–Alegre [5] for the line-digraph method in the $(d, D)$ problem; see also the modern overview [6]. This theorem along with our recursive formula for $\rho$ allows us to calculate both $\rho$ and $\sigma$ recursively.

**Theorem 1.1.** *Let $G = K(d, D)$ and $L(G) \cong K(d, D+1)$ be its line digraph. Then for each $1 \leq k \leq D$,*
$$\sigma_{k+1}(L(G)) = d^2 \sigma_k(G) - \rho_k(G).$$

*Proof.* Let $E(G)$ be the directed edges of $G$. Vertices of $L(G)$ are the directed edges of $G$. Write an edge of $G$ as $(u \to v)$, where $u$ and $v$ are length $(D+1)$ Kautz words with $v$ obtained from $u$ by shifting left one position and appending a new letter. Then there is an edge $e \to f$ in $L(G)$ precisely when the second word of $e$ equals the first word of $f$; that is, if $e = (u \to v)$ and $f = (v \to w)$. Equivalently, if $u = a_0 a_1 \cdots a_D$, then $v = a_1 a_2 \cdots a_{D+1}$ and $w = a_2 a_3 \cdots a_{D+2}$, so $e \to f$ corresponds to the overlap

$$(a_0 a_1 \cdots a_D \to a_1 a_2 \cdots a_{D+1}) \to (a_1 a_2 \cdots a_{D+1} \to a_2 a_3 \cdots a_{D+2}).$$

A path $e = e_0, \ldots, e_{k+1} = f$ of length $k+1$ in $L(G)$ corresponds to a path of length $k+2$ in $G$, so $\text{dist}_{L(G)}(e, f) = k+1$ if and only if $\text{dist}_G(u_1, u_{k+1}) = k$ for $e = (u_0, u_1)$, $f = (u_{k+1}, u_{k+2})$. For each $(u, v)$ at distance $k$ there are $d^2$ pairs $(e, f)$ with $e = (x, u)$, $f = (v, y)$. Exactly one pair *completes* the shortest directed cycle through the middle edge at length $k+1$, removing one from the count; summing over all $(u, v)$ yields $d^2 \sigma_k - \rho_k = \sigma_{k+1}(L(G))$. □

**Lemma 1.2.** *Fix $D \geq 1$ and $d \geq 2$. Then*
$$\Delta_k(d, D) = 0 \quad \text{for } k < \lfloor D/2 \rfloor + 2,$$
*and $\Delta_k(d, D)$ is non-negative for every $k$ with $\lfloor D/2 \rfloor + 2 \leq k \leq D - 1$.*

*Proof.* Write vertices of $K(d, D)$ as Kautz words $w_0 \cdots w_D$ over $[q] = \{0, \ldots, d\}$ with $w_i \neq w_{i+1}$. For an oriented edge $e = (u \to v)$, let
$$r := \max\Big\{ t \in \{0, \ldots, D+1\} : \text{the suffix of } v \text{ of length } t \text{ equals the prefix of } u \text{ of length } t\Big\}.$$
Then the shortest return path $v \to u$ has length $D+1-r$, so
$$|Z_e| = 1 + (D+1-r) = D+2-r, \quad \text{hence} \quad k = |Z_e| - 1 = D+1-r,$$
and we call $r$ the overlap size.

*(1) Persistence from $D-1$ to $D$.* Form $K(d, D) \cong L(K(d, D-1))$, the line digraph: vertices are directed edges of $K(d, D-1)$ and adjacency encodes head–tail matching of those edges. A cycle in $K(d, D-1)$ containing an edge $e$ maps to a cycle of the *same length* through the corresponding vertex of $L(K(d, D-1)) \cong K(d, D)$, and any cycle in the line digraph that projects through that vertex maps back to a cycle in $K(d, D-1)$ of the same length. Hence a shortest cycle through $e$ at diameter $D-1$ remains shortest at diameter $D$; no shorter cycle can suddenly appear at $D$ (its projection would contradict minimality at $D-1$). Therefore every shortest cycle at diameter $D-1$ persists at $D$.

*(2) Vanishing window $k < \lfloor D/2 \rfloor + 2$.* New first hits at diameter $D$ can arise only when the maximal suffix/prefix overlap between the endpoint words first *uses* the new rightmost symbol introduced when passing from $D-1$ to $D$. Let $r$ be that (first) overlap length at row $D$, and recall the edge-length relation
$$|Z_e| = (D+2) - r, \quad \text{so} \quad k = |Z_e| - 1 = (D+1) - r.$$
If $k \leq \lfloor D/2 \rfloor + 1$, then $r \geq D+1 - (\lfloor D/2 \rfloor + 1) = \lceil D/2 \rceil$. But then already at row $D-1$ the two words shared an overlap of size at least $r - 1 \geq \lceil D/2 \rceil - 1$, which yields a return of length
$$(D+1) - (r-1) = D+2-r = k+1.$$
Thus the same (or shorter) cycle would have existed at diameter $D-1$, so it cannot be *new* at $D$. Consequently, no new first hits occur in row $D$ for $k \leq \lfloor D/2 \rfloor + 1$, i.e.,
$$\Delta_k(d, D) = 0 \quad \text{for all } k < \lfloor D/2 \rfloor + 2.$$

This proves persistence and the stated vanishing range for all $d \geq 2$. □

We shall show that $\Delta_k(d, D)$ is positive for every $k$ with $\lfloor D/2 \rfloor + 2 \leq k \leq D - 1$ and give a recursion which calculates it exactly.

## 2 Terminal letter state space

Fix $d \geq 2$ and set $q = d+1$. The *terminal letter state space* is

$$\mathcal{S}_1 = \{[\ell, r] \in [q]^2 : \ell \neq r\}, \qquad n := |\mathcal{S}_1| = q(q-1).$$

We index $\mathcal{S}_1$ in lexicographic order and write $\mathbf{1} \in \mathbb{N}^n$ for the all-ones column. A single legal Kautz extension appends one symbol on the right and shifts the terminal pair; this is the *one-letter transfer matrix*

$$S_q\big([\ell, r] \to [x, y]\big) := \mathbf{1}\{x = r,\ y \neq r\},$$

so each one-letter advance moves $[\ell, r] \mapsto [r, y]$ with exactly $d = q-1$ choices for $y$. This matrix is identical to the adjacency matrix of $K(d, 2)$.

To enforce "*first hit at exactly $k$*" within the same $D$-row, we interleave one-letter transfers with small diagonal *boundary masks* $\Lambda^{(s)}$ that (i) forbid steps that would realize a shorter return while the current offset to the new right boundary is $s > 0$, and (ii) *force* closure at the last step via $(I - \Lambda^{(0)})$. Far from the boundary the evolution is unconstrained, so $\Lambda^{(s)} = I$ for all sufficiently large $s$ (in particular, for $d = 2$ only $s \in \{0, 1, 2\}$ are nontrivial).

Let a unit row $e^\top$ select any start state in $\mathcal{S}_1$, and define the offsets

$$s_j := (D - 1) - k + j \qquad (j = 1, \ldots, k-1),$$

i.e., the distance to the boundary after $j-1$ one-letter advances.

**Theorem 2.1.** *For any $d \geq 2$, $D \geq 1$, and $1 \leq k \leq D$,*

$$\Delta_k(d, D) = n \cdot e^\top \Big( \prod_{j=1}^{k-1} S_q \Lambda^{(s_j)} \Big) S_q (I - \Lambda^{(0)}) \mathbf{1}, \quad \text{with } n = q(q-1).$$

*Here $S_q$ is the transfer matrix on $\mathcal{S}_1$, the diagonal masks $\Lambda^{(s)}$ depend only on local checks at offset $s$, and $\Lambda^{(s)} = I$ for all sufficiently large $s$.*

*Proof.* Consider an oriented edge $e = (u \to v)$ in $K(d, D)$, written as a terminal pair $[\ell, r]$ of the word $u$ at the right end. A one-letter update adds one symbol on the right of $v$ while the leftmost symbol of $u$ falls off; this is exactly the linear action of $S_q$. If after $j-1$ such updates the offset to the new right boundary is $s_j > 0$, then any step that would complete a strictly shorter return must be excluded; this is enforced by the diagonal mask $\Lambda^{(s_j)}$. After $k-1$ constrained one-letter steps, the last step must *hit* the boundary (i.e., close the shortest return at length $k+1$), which is enforced by $I - \Lambda^{(0)}$. Thus, for a fixed start row $e^\top$, the scalar

$$e^\top \Big( \prod_{j=1}^{k-1} S_q \Lambda^{(s_j)} \Big) S_q (I - \Lambda^{(0)}) \mathbf{1}$$

counts exactly the terminal-pair paths realizing a first return of length $k+1$. By vertex transitivity of $K(d, 2)$, the per-start count is independent of the chosen start state so summing over all $n = q(q-1)$ starts multiplies by $n$, proving the formula. □

**Value of $\rho_D(d, D)$**

The counts $\{\rho_k(d, D)\}_{k=1}^{D}$ partition the edges:

$$|E(K(d, D))| = (d+1)\, d^{D+1} \;=\; \sum_{k=1}^{D} \rho_k(d, D).$$

. Thus

$$\rho_D(d, D) \;=\; (d+1)\, d^{D+1} \;-\; \sum_{k=1}^{D-1} \bigl[\rho_k(d, D-1) + \Delta_k(d, D)\bigr].$$

## 3   Specialization to $d = 2$: explicit matrices, masks, and examples

For $d = 2$ set $q = 3$. The *terminal letter state space* is

$$\mathcal{S}_1 = \{[\ell, r] \in [q]^2 : \ell \neq r\}, \qquad n := |\mathcal{S}_1| = q(q-1) = 6.$$

We index $\mathcal{S}_1$ in the canonical order $[01], [02], [10], [12], [20], [21]$, and write $\mathbf{1} \in \mathbb{N}^6$ for the all-ones column vector.

**One-letter transfer matrix $d = 2$.**  In the above state order,

$$S_3 \;=\; \begin{bmatrix} 0 & 0 & 1 & 1 & 0 & 0 \\ 0 & 0 & 0 & 0 & 1 & 1 \\ 1 & 1 & 0 & 0 & 0 & 0 \\ 0 & 0 & 0 & 0 & 1 & 1 \\ 0 & 0 & 1 & 1 & 0 & 0 \\ 1 & 1 & 0 & 0 & 0 & 0 \end{bmatrix}, \qquad S_3\,\mathbf{1} \;=\; 2\,\mathbf{1}, \quad \mathbf{1}^\top S_3 \;=\; 2\,\mathbf{1}^\top.$$

Thus every interior one-step contributes a factor 2 in row/column sums.

**Boundary masks for $d = 2$.**  Let $\Lambda^{(s)}$ be diagonal masks encoding "no early hit" at offset $s > 0$ from the new right boundary, and let $(I - \Lambda^{(0)})$ enforce "hit now" at the final step. For $d = 2$ only $s \in \{0, 1, 2\}$ are nontrivial; in the canonical order:

$$\Lambda^{(0)} = \mathrm{diag}(0,0,1,0,1,0), \quad \Lambda^{(1)} = \mathrm{diag}(1,0,1,1,0,1), \quad \Lambda^{(2)} = \mathrm{diag}(1,1,0,1,0,1), \quad \Lambda^{(s)} = I\ (s \geq 3).$$

These masks depend only on local checks near the boundary; away from the boundary the evolution is unconstrained.

**Calculation of $\Delta_k(2, D)$.**  With $d = 2$, $D \geq 1$, and $1 \leq k \leq D$ then

$$\Delta_k(2, D) \;=\; 6 \,\cdot\, e^\top \Bigl(\prod_{j=1}^{k-1} S_3\, \Lambda^{(s_j)}\Bigr) S_3\, (I - \Lambda^{(0)})\, \mathbf{1}, \qquad s_j = D - k + j.$$

**Examples.** Write $r^* = D - k$ (initial offset to the boundary); the schedule is $k$ single-letter steps with masks at offsets $s_j = D - k + j$.

- **Row** $D = 10$ (active $k = 7, 8, 9$): $k = 9$: offsets $s = 1, 2, 3, 4, 5, 6, 7, 8$ before the final hit; only $s \in \{1, 2\}$ are masked. Per start the product gives 27, hence $\Delta_9(2, 10) = 6 \cdot 27 = \mathbf{162}$. $k = 8$: masked at $s = 2$ then interior, per start $6 \Rightarrow \Delta_8(2, 10) = \mathbf{36}$. $k = 7$: two near-boundary masks, per start $2 \Rightarrow \Delta_7(2, 10) = \mathbf{12}$.

- **Row** $D = 11$ (active $k = 7, 8, 9, 10$): $k = 10$: four unconstrained near-boundary steps give per start $64 \Rightarrow \Delta_{10}(2, 11) = \mathbf{384}$. $k = 9$: per start $16 \Rightarrow \Delta_9(2, 11) = \mathbf{96}$. $k = 8$: per start $3 \Rightarrow \Delta_8(2, 11) = \mathbf{18}$. $k = 7$: forced, per start $1 \Rightarrow \Delta_7(2, 11) = \mathbf{6}$.

- **Far right entry** $D = 15$, $k = 12$: Per start $60 \Rightarrow \Delta_{12}(2, 15) = 6 \cdot 60 = \mathbf{360}$.

The masks affect only the last few steps near $s = 2, 1, 0$; all prior steps are pure $S_3$.

## Enumerating $\rho_D(d, D)$ from primitive necklaces

Set $q = d+1$ and $n = D+1$. A directed cycle of length $n$ in $K(d, D)$ corresponds to a rotational equivalence class of proper $q$-colorings of the cycle graph $C_n$ (equivalently, cyclic words over $[q]$ with no equal adjacent symbols, including the wrap). Let $\mathcal{N}_{\text{prim}}(n; q)$ denote the number of *primitive* (aperiodic) such colorings. By Möbius inversion,

$$\mathcal{N}_{\text{prim}}(n; q) = \frac{1}{n} \sum_{j | n} \mu(j) (q-1)^{n/j}.$$

Each primitive $n$-cycle contributes exactly $n$ oriented edges whose shortest directed cycle has length $n$. Therefore $\rho_D(d, D)$ (i.e., shortest cycle length $D+1 = n$) is

$$\rho_D(d, D) = n \cdot \mathcal{N}_{\text{prim}}(n; q) = \sum_{j | (D+1)} \mu(j) (q-1)^{\frac{D+1}{j}}.$$

This identity applies only to $\rho_D$; we are not aware if $\rho_k$ with $k < D$ has ever been enumerated. See Lyndon [3] and Smyth [4] for classical background on aperiodic (primitive) necklaces and Möbius inversion.

## References

[1] W. H. Kautz. Bounds on Directed (d,k) Graphs. *Theory of Switching, Part 1*, 1968.

[2] M. Aigner. On the linegraph of a directed graph. *Mathematische Zeitschrift*, 102:56–61, 1967.

[3] R. C. Lyndon. On Burnside's problem. *Transactions of the American Mathematical Society*, 77(2):202–215, 1954.

[4] W. F. Smyth. *Computing Patterns in Strings*. Pearson Addison–Wesley, 2003.

[5] M. A. Fiol, J. L. A. Yebra, and I. Alegre. Line digraph iterations and the $(d, k)$ problem for directed graphs. *(Proc./journal version mid-1980s; insert venue/year)*.

[6] C. Dalfó, M. A. Fiol, M. Miller, J. Ryan, and J. Širáň. An algebraic approach to lifts of digraphs. *arXiv:1612.08855*, 2016.